\documentclass[12pt]{article} \textwidth=6in \oddsidemargin=0in
\usepackage{amssymb}

\begin{document}

\def\bc{\begin{center}} \def\ec{\end{center}}
\def\ov{\over} \def\ep{\varepsilon} 
\newcommand{\C}[1]{{\cal C}_{#1}} \def\inv{^{-1}}\def\be{\begin{equation}} \def\noi{\noindent}
\def\ee{\end{equation}}\def\x{\xi}\def\({\left(} \def\){\right)} \def\iy{\infty} \def\e{\eta} \def\cd{\cdots} \def\ph{\varphi} \def\ps{\psi} \def\La{\Lambda} \def\s{\sigma} \def\ds{\displaystyle}
\newcommand{\xs}[1]{\x_{\s(#1)}} \newcommand\sg[1]{{\rm sgn}\,#1}
\newcommand{\xii}[1]{\x_{i_{#1}}} \def\ld{\ldots} \def\Z{\mathbb Z}
\def\a{\alpha}\def\g{\gamma}\def\b{\beta} \def\noi{\noindent}
\def\ar{\leftrightarrow} \def\d{\delta} \def\S{\mathbb S} \def\tl{\tilde} \def\E{\mathbb E} \def\P{\mathbb P} \def\TS{U\backslash T}\def\CS{{\cal S}} \newcommand{\br}[2]{\left[{#1\atop #2}\right]}
\def\bs{\backslash} \def\CD{{\cal D}} \def\N{\mathbb N} \def\bC{\mathbb{C}} \newcommand{\vs}[1]{\vskip#1ex} \def\t{\tau} \def\l{\ell}
\def\la{\lambda} \def\ch{\raisebox{.3ex}{$\chi$}} \def\phy{\ph_\iy} \def\m{\mu} \def\KA{K_{\rm Airy}} \def\A{{\rm{Ai}}} \def\z{\zeta}
\def\tr{\textrm{tr}} \def\G{\Gamma} \def\ps{\psi} \def\dl{\delta}
\def\k{\kappa} \def\c{\cal{c}} \def\r{\rho} \def\phy{\ph_\iy}
\def\si{\s\inv}

\bc {\Large\bf Formulas and Asymptotics for\\
\vs{1}
the Asymmetric Simple Exclusion Process}\footnote{This is an expanded version of a series of lectures delivered by the second author at Universit\'e de Paris in June, 2009, describing the results in the articles \cite{TW1,TW2,TW3}. Although complete proofs will in general not be presented here, at least the main elements of them will be.}  \ec

\vs{2}
\begin{center}{\large\bf Craig A.~Tracy}\\
{\it Department of Mathematics \\
University of California\\
Davis, CA 95616, USA}\end{center}

\begin{center}{\large \bf Harold Widom}\\
{\it Department of Mathematics\\
University of California\\
Santa Cruz, CA 95064, USA}\end{center}

\bc{\bf I. Introduction}\ec

The asymmetric simple exclusion process (ASEP) is a special case of processes introduced in 1970 by F. Spitzer \cite{Sp}. In ASEP, particles are at integer sites on the line. Each particle waits exponential time, and then
\vs{.5}
(1) with probability $p$ it moves one step to the right if the site is unoccupied, otherwise it does not move;
\vs{.5}
(2) with probability $q=1-p$ it moves one step to the left if the site is unoccupied, otherwise it does not move.

In the totally asymmetric simple exclusion process (TASEP) particles can only move in one direction, so either $p=0$ or $q=0$. In a major breakthrough, K.~Johansson \cite{J} related a probability in TASEP to a probability in random matrix theory. If $q=0$ and initially particles were at the negative integers, then the probability that at time~$t$  the particle initially at $-m$ has moved at least  
$n\ (\ge m)$ times equals the probability distribution for the largest eigenvalue in the Laguerre ensemble of $m\times m$ matrices with weight function $x^{n-m}\,e^{-x}$. Thus, it is given by a constant depending on $m$ and $n$ times the determinant
\[\det\(\int_0^t x^{n-m+i+j}\,e^{-x}\,dx\)_{i,\,j=0,\ldots,m-1}.\]
This connection led to considerable progress in understanding TASEP and the derivation of asymptotic results. For ASEP there is no longer a determinantal structure and a different approach was required. 

Here are the main results of \cite{TW1,TW2,TW3}.
First we consider ASEP with finitely many particles. For $N$-particle ASEP a possible configuraion is given by 
\[X=\{x_1,\ld,x_N\},\quad x_1<\cd<x_N,\quad (x_i\in\Z).\]
Thus the $x_i$ are the occupied sites. We assume an initial configuration $Y=\{y_1,\ld,y_N\},$
and obtain formulas for
\vs{.5}
(1) $\P_Y(X;t)$, the probability that at time $t$ the system is in configuration~$X$.\footnote{This had been known for the case $N=2$ \cite{Sc}.}
\vs{.5}
(2) $\P_Y(x_m(t)=x)$, the probability that at time $t$ the $m$th particle from the left is at $x$.  
\vs{.5}
The formula we get for the latter  extends to infinite systems 
\[y_1<y_2<\cd\to +\iy.\]
In particular we may take $Y=\Z^+$. (This is the {\it step initial condition.})

For the derivation of (1) we use the {\it Bethe Ansatz} \cite{B} to obtain a solution of a differential equation with boundary conditions. The derivation of (2) from (1) requires the proof of two combinatorial identities. The derivation we outline for this is from \cite{TW4} and simpler than the one in \cite{TW1}.

For step initial condition we show that $\P(x_m(t)=x)$ has a representation in terms of Fredholm determinants. This makes asymptotic analysis possible. We assume that $q>p$, so there is a drift to the left, define $\g=q-p$, and obtain asymptotics as $t\to\iy$ for
\[\P\left(x_m(t/\g)>x\right)\ \ \ (\textrm{fixed $m$ and $x$)},\]
and the limits as $t\to\iy$ of
\[\P\(x_m(t/\g)\le -t+\g^{1/2}\,s\,t^{1/2}\)\ \ \ (\textrm{$m$ fixed)},\]
\[\P\(x_m(t/\g )\le-c_1(\s) \,t+c_2(\s)\,s\,t^{1/3}\)\ \ \ (m=[\s t]),\]
where $c_1(\s)$ and $c_2(\s)$ are certain explicit constants. 

The last limit is the distribution function $F_2(s)$ of random matrix theory, the limiting distribution for the rescaled largest eigenvale in the Gaussian unitary ensemble. These asymptotics  were obtained in \cite{J} for the case of TASEP. (That $F_2$ should arise in ASEP had long been suspected. This is referred to as 
\textit{KPZ universality} \cite{KPZ}.)

\pagebreak

\bc{\bf II. Integral Formulas}\ec

\bc{\bf 1. The differential equation}\ec

The idea goes back to \cite{B}. There is a differential equation with boundary conditions and an initial condition whose solution gives $P_Y(X;t)$. To state it we introduce the new notation $u(X;t)$ or $u(X)$ in place of $P_Y(X;t)$.\footnote{The reason is that if $X=\{x_1,\ld,x_N\}$ then $P_Y(X;t)$ only makes sense when $x_1<\cd<x_N$, but for  $u(X;t)$ there will be no such requirement.}

We first consider the case $N=2$, and consider $du(x_1,x_2)/dt$. After an exponential waiting time, the system could enter state $\{x_1,\,x_2\}$ or it could leave this state. Assume first that $x_2-x_1>1$, so that there is no interference between the two particles. The system could enter the state if the first particle had been at $x_1-1$ (this has probability $u(x_1-1,\,x_2)$) and moved one step to the right (probability $p$), and three other analogous ways. The system could leave the state if the first particle is at $x_1$ (probability $u(x_1,\,x_2)$) and moves one step to the right (probability $p$) or one step to the left (probability $q$), and analogously for the second particle. These give the equation

\[{d\ov dt}\,u(x_1,x_2)=p\,u(x_1-1,x_2)+q\,u(x_1+1,x_2)\]
\[+ p\,u(x_1,x_2-1)+q\,u(x_1,x_2+1)-2\,u(x_1,x_2).\]

But if $x_2-x_1=1$ then for entering the state the first particle could not have been one step to the right nor the second particle one step to the left, and for leaving the state the first particle cannot move right nor can the second particle move left. Therefore in this case
\[{d\ov dt}\,u(x_1,x_2)=p\,u(x_1-1,x_2)+q\,u(x_1,x_2+1)-
u(x_1,x_2).\]

We could combine these two equations into one, but then the right side would have nonconstant coefficients. Instead, as in \cite{B}, we observe that if we formally subtract the two equations we get, when $x_2=x_1+1$,
\[0=p\,u(x_1,x_1)+q\,u(x_1+1,x_1+1)-u(x_1,x_1+1).\]
If the first equation holds for {\it all} $x_1$ and $x_2$, and this last {\it boundary condition} holds for all $x_1$, then the second equation holds when $x_2=x_1+1$. So an equation with nonconstant coefficients has been replaced with an equation with constant coefficients plus a boundary condition.

This was done for $N=2$, but it holds for general $N$. Suppose the function $u(X;t)$, defined for all 
$X=\{x_1,\ld,x_N\}\in\Z^N$, satisfies the {\it master equation}

\[{d\ov dt}\,u(X;t)\]
\[=\sum_{i=1}^N\,[p\,u(\ld,x_i-1,\ld)+q\,u(\ld ,x_i+1,\ld)-u(X)],\]
and the {\it boundary conditions} 
\[u(\ld,x_i,x_i+1,\ld)=p\,u(\ld,x_i,x_i,\ld)+q\,u(\ld,x_i+1,x_i+1,\ld) \]
for $i=1,\ld,N-1$.\footnote{For $N\ge3$ the boundary conditions arising from configurations with more than two adjacent particles automatically follow from the boundary conditions arising from two adjacent particles.} Suppose also that it satisfies the {\it initial condition}
\[u(X;0)=\dl_Y(X)\ \ {\rm when}\ x_1<\cd<x_N,\]
which reflects the initial configuration $Y$. Then 
\[u(X;t)=P_Y(X;t)\ \ {\rm when}\ x_1<\cd<x_N.\]

Thus the strategy will be: (1) find a large class of solutions to the master equation; (2) find a subset satisfying the boundary conditions; (3) find one of these satisfying the initial condition. The last will be the hard (and new) part.

\bc{\bf 2. Solutions to the master equation}\ec

Define
\[\ep(\x)=p\,\x\inv+q\,\x-1.\]
For any nonzero complex numbers $\x_1,\ld,\x_N$, a solution of the equation is
\[\prod_i\(\x_i^{x_i}\,e^{\ep(\x_i)\,t}\).\]
Since the $\x_i$ are arbitrary another solution is obtained by permuting them. Thus, for any $\s$ in the symmetric group $\S_N$  another solution is
\[\prod_i\x_{\s(i)}^{x_i}\,\prod_i e^{\ep(\x_i)\,t}.\]
(The second factor is symmetric in the $\x_i$, which is why we can write it as we do.) Since the equation is linear, any linear combination of these is a solution, as is any integral (over the $\x_i$) of a linear combination. Thus we arrive at the {\it Bethe Ansatz solutions} 
\be u(X;t)=\int \sum_{\s\in \S_N} F_\s(\x)\,\prod_i\xs{i}^{x_i}\,\prod_i e^{\ep(\x_i)\,t}\,d^N\x.\label{u}\ee
The $F_\s$ are arbitrary functions of the $\x_i$, and the domain of integration is arbitrary.

\bc{\bf 3. Satisfying the boundary conditions}\ec 

We look for functions $F_\s$ such that the integrand satisfies the boundary conditions pointwise. The $i$th boundary condition is satisfied pointwise when 

\be\sum_{\s\in\S_N}\; F_\s\,(p+q\,\xs{i}\,\xs{i+1}-\xs{i+1})\; (\xs{i}\xs{i+1})^{x_i}\,\prod_{j\ne i,\,i+1}\xs{j}^{x_j}=0.\label{bdrycond}\ee

Define $T_i\s$ to be the permutation that differs from $\s$ by an intechange of the $i$th and $(i+1)$st entries. Thus, if $\s=(2\ 3\ 1\ 4)$ then $T_2\s=(2\ 1\ 3\ 4)$. Since $T_i$ is bijective, (\ref{bdrycond}) is unchanged if each $\s$ in the summand is replaced by $T_i\s$, and therefore unchanged if we add the two. Since the last two factors are unchanged upon replacing $\s$ by $T_i\s$, we see that a sufficient condition that (\ref{bdrycond}) satisfied is that for each $\s$,
\[F_\s\,(p+q\xs{i}\xs{i+1}-\xs{i+1})+F_{T_i\s}\,(p+q\xs{i}\xs{i+1}-\xs{i})=0.\]
Because the expression will appear so often it is convenient to define
\[f(\x,\,\x')=p+q\x\x'-\x,\]
so the condition becomes
\[F_\s\,f(\xs{i+1},\,\xs{i})+F_{T_i\s}\,f(\xs{i},\,\xs{i+1})=0.\]

This is to hold for all $\s$ and all $i$. Since these are $(n-1)\,n!$ linear equations in the $n!$ unknowns $F_\s$, we cannot necessarily expect a solution. But there are solutions, and in fact it is easy to see that 
\[F_\s(\x)=\sg\s\,\prod_{i<j}f(\xs{i},\,\xs{j})\times \ph(\x),\]
where $\ph$ is an arbitrary function of the $\x_i$, satisfies the equations. (In fact this is the general solution, since if any $F_\s$ is known then all others are determined.) These $F_\s$ in (\ref{u}) give a family of solutions to the master equation that satisfy the boundary conditions. 
\pagebreak

\bc{\bf 4. Satisfying the initial condition}\ec

The initial condition is
\be\int \sum_{\s\in \S_N} F_\s(\x)\,\prod_i\xs{i}^{x_i}\,d^N\x=\dl_Y(X).\label{initcond}\ee
We begin with the fact that if $\C{}$ is a contour enclosing zero then\footnote{Unless specified otherwise, all contours are described once counterclockwise, and all contour integrals have a factor $1/2\pi i$.}
\[\int_{\C{}} \x^{x-y-1}\,d\x=\dl_y(x).\]
Therefore
\[\int_{\C{}^N}\prod_i \x_i^{x_i-y_i-1}\,d^N\x=\dl_Y(X).\]
Thus if $id$ denotes the identity permutation, then the $\s=id$ summand in (\ref{initcond}) will give the integral $\dl_Y(X)$ if the integration is over $\C{}^N$ and
\[F_{id}(\x)=\prod_i \x_i^{-y_i-1}.\]

For this to hold we choose 
\[\ph(\x)=\prod_{i<j}f(\x_{i},\,\x_{j})\inv\cdot\prod_i \x_i^{-y_i-1}.\] 
If we define
\be A_\s=\sg\s\,{\prod_{i<j}f(\xs{i},\,\xs{j})\ov \prod_{i<j}f(\x_i,\,\x_j)}\label{A}\ee
then the solution we have chosen is
\be u(X;t)=\sum_\s \int_{\C{}^N} A_\s(\x)\,\prod_i\xs{i}^{x_i}\,\prod_i \Big(\x_i^{-y_i-1}\,e^{\ep(\x_i)\,t})\,d^N\x.\label{sol}\ee
It satisfies the master equation and boundary conditions, and the $\s=id$ summand satisfies initial condition. 

Observe that because of the poles of $A_\s$ when $\s\ne id$ this will not be well-defined until we specify $\C{}$ further.

\bc{\it a. TASEP}\ec

When $p=1$ we have
\[A_\s={\rm sgn}\,\s\prod_i(1-\x_{\s(i)})^{\s(i)-i}.\]
Because of this product structure the integrals of products in (\ref{sol}) may be written as product of integrals and (\ref{sol}) becomes
\pagebreak
\[u(X;t)=\sum_\s\sg\s\,\prod_i\int_{\C{}}(1-\x^{\s(i)-i})\,\x^{x_i-y_{\s(i)}-1}\,e^{(\x\inv-1)t}\,d\x\]
\[=
\det\(\int_{\C{}}(1-\x)^{j-i}\,\x^{x_i-y_{j}-1}e^{(\x\inv-1)t}\,d\x\).\]
Sch\"utz \cite{Sc} obtained this solution to the master equation, using Bethe Ansatz as we have described, and went further to show that it satisfies the boundary condition when the point $\x=1$ is outside the contour $\C{}$. So he established the formula
\[P_Y(X;t)=\det\(\int_{\C{r}}(1-\x)^{j-i}\,\x^{x_i-y_{j}-1}e^{(\x\inv-1)t}\,d\x\),\]
where $\C{r}$ denotes the circle with center zero and radius $r<1$.

\begin{center}{\it b. ASEP}\end{center} 

In \cite{Sc} Sch\"utz also considered ASEP and showed that when $N=2$ the probability $P_Y(X;t)$ is equal to a sum of a two-dimensional integral and a one-dimensional integral. In the two-dimensional integral the contours were different. It turns out that if one integrates over small contours only then the sum is the sum of two two-dimensional integrals. And this extends to general $N$.

Recall that because of the poles of $A_\s$, it matters which contours $\C{}$ we take in  (\ref{sol}). When $p\ne0$ all poles of the $A_\s$ will lie outside $\C{r}$ if $r$ is small enough. These are the contours we take. 
\vs{2}
{\bf Theorem} \cite[Theorem~2.1]{TW1}: Suppose $p\ne0$ and assume that $r$ is so small that all poles of the $A_\s$ lie outside $\C{r}$. Then
\be P_Y(X;t)=\sum_\s \int_{\C{r}^N} A_\s(\x)\,\prod_i\xs{i}^{x_i}\,\prod_i \Big(\x_i^{-y_i-1}\,e^{\ep(\x_i)\,t})\,d^N\x.\label{PYX}\ee

For the proof we have to show that the initial condition is satisfied. Since the $\s=id$ summand satisfies the initial condition, what is to be shown is that if
\[I(\s)=\int_{\C{r}^N} A_\s(\x)\,\prod_i\xs{i}^{x_i}\,\prod_i \x_i^{-y_i-1}\,d^N\x,\]
then 
\be\sum_{\s\ne id} I(\s)=0.\footnote{There was an error in the proof of this in \cite{TW1}. A corrected version of \cite{TW1} is posted as arXiv:0704.2633. We outline the proof here.}\label{Isum}\ee
 
We show that the permutations in $\S_N\bs\{id\}$ can be grouped in such a way that the sum of the $I(\s)$ fron each group is equal to zero. For $1\le n<N$ fix $n-1$ distinct numbers $i_1,\ld,i_{n-1}\in[1,\,N-1]$, define $A=\{i_1,\ld,i_{n-1}\}$,
and then
\[\S_N(A)=\{\s\in\S_N:\s(1)=i_1,\ld,\s(n-1)=i_{n-1},\,\s(n)=N\}.\]
When $n=1$ these are all permutations with $\s(1)=N$. When $n=N-1$ each $\S_N(A)$ consists of a single permutation. Let $B$ be the complement of 
$A\cup\{N\}$  in $[1,\,N]$. 

We first we make the substitution
\[\x_N\to{\e\ov\prod_{i<N}\x_i}\]
in all the integrals. The product of the powers of the $\x_i$ in (\ref{PYX}) becomes
\be\e^{x_n-y_N-1}\,\prod_{i<N}\x_i^{x_{\s\inv(i)}-x_n+y_N-y_i-1}.\label{xiprod}\ee
We use the alternative representation
\be A_\s=\prod_{\l<k\atop\si(\l)>\si(k)}\!\!\left(-{f(\x_k,\,\x_\l)\ov f(\x_\l,\,\x_k)}\right)\label{AltA}\ee
to see what happens when we shrink some of the $\x_i$-contours. The only poles we might cross come from the denominatiors in (\ref{AltA}) after the substitution, and these are in the $\x_i$-variables when $i\in B$.
\vs{1}
\noi{\bf Lemma 1}. When $n=N-1$ we have $I(\s)=0$ for $\s\in\S_N(A)$.
\vs{1}
There is a single $i\in B$ and in this case there is no pole in the $\x_i$-variable coming from (\ref{AltA}). Using $x_N>x_{N-1}$ and $y_N>y_{i}$ in (\ref{xiprod}), we see that the integrand is analytic at $\x_i=0$. Therefore the integral with respect to $\x_i$ is zero.
\vs{1}
\noi{\bf Lemma 2}. When $n<N-1$ all $I(\s)$ with $\s\in\S_N(A)$ are sums of lower-order integrals in each of which a partial product in (\ref{AltA}) independent of $\s\in A$ is replaced by another factor. In each integral some $\x_i$ with $i\in B$ is equal to another $\x_j$ with $j\in B$.
\vs{1}

If $j=\max B$, we shrink the $\x_j$-contour and obtain $(N-1)$-dimensional integrals coming from poles associated with the variables $\x_k$ with $k\in B\bs\{j\}$. For each such $k$ we integrate with respect to $\x_k$ the residue at this pole by shrinking the contour, and obtain $(N-2)$-dimensional integrals having the property described in the lemma.
\vs{1}
\noi{\bf Lemma 3}. For each integral of Lemma 2 there is a partition of $\S_N(A)$ into pairs $\s,\,\s'$ such that $I(\s)+I(\s')=0$ for each pair.
\vs{1}

Consider an integral in which $\x_i=\x_j$. We pair $\s$ and $\s'$ if $\s\inv(i)={\s'}\inv(j)$ and $\s\inv(j)={\s'}\inv(i)$, and $\s\inv(k)={\s'}\inv(k)$ when $k\ne i,j$. The factor (\ref{xiprod}) is clearly the same for both when $\x_i=\x_j$, and $A_\s$ and $A_{\s'}$ are negatives of each other then.

Here is why. Assume for definiteness that
$i<j$ and $\si(i)<\si(j).$ Then the factor 
corresponding to $\l=i$ and $k=j$ does not appear for $\s$ in (\ref{AltA}) but it does appear for $\s'$.  This factor equals $-1$ when $\x_i=\x_j$. And it is straightforward to check that for any $k\ne i,j$ the product of factors involving $k$ and either $i$ or $j$ is the same for $\s$ and $\s'$ when $\x_i=\x_j$.

\vs{1}
Now (\ref{Isum}) can be shown by induction on $N$. When $N=2$ it folows from Lemma~1. Assume $N>2$ and that the result holds for $N-1$. For those permutations for which $\s(N)=N$ we integrate with respect to $\x_1,\ld,\x_{N-1}$ and use the induction hypothesis. The set of permutations with $\s(N)<N$ is the disjoint union of the various $\S_N(A)$, and for these we apply Lemmas 1 and 3.

\begin{center}{\bf 5. The left-most particle}\end{center} 

The probability $\P_Y(x_1(t)=x)$ is the sum of $\P(X;t)$ over all $X$ for which $x_1=x$, thus over all $x_2,\ld,x_N$ statisfying $x<x_2<\cd <x_N<\iy.$ When $r<1$ we may sum under the integral sign in (\ref{PYX}), and the integrand becomes

\[{\prod_i(\x_i^{x-y_i-1}e^{\ep(\x_i)t})\ov\prod_{i<j}f(\x_i,\,\x_j)}
\cdot \sum_{\s}\;\sg\s\,\(\prod_{i<j}f(\x_{\s(i)},\,\x_{\s(j)})\right.\]
\[\left.\times\,
{\x_{\s(2)}\x_{\s(3)}^2\cd\x_{\s(N)}^{N-1}\ov (1-\x_{\s(2)}\x_{\s(3)}\cd\x_{\s(N)}) (1-\x_{\s(3)}\cd\x_{\s(N)})\cd(1-\x_{\s(N)})}\).\]

Fortunately we have our first combinatorial identity,\footnote{Doron Zeilberger saw the identity when it was still a conjecture and suggested to the authors that an identity of I.~Schur  \cite[Problem VII.47]{PSz} had a similar look about it and might be proved in a similar way. This led to the proof we present.}

\[ \sum_{\s}\,\sg\s \(\prod_{i<j}f(\x_{\s(i)},\,\x_{\s(j)})\cdot 
{\x_{\s(2)}\x_{\s(3)}^2\cd\x_{\s(N)}^{N-1}\ov (1-\x_{\s(1)}\x_{\s(2)}\cd\x_{\s(N)}) (1-\x_{\s(2)}\cd\x_{\s(N)})\cd(1-\x_{\s(N)})}\)\]
\be=p^{N(N-1)/2}\, {\prod_{i<j} (\x_j-\x_i)\ov \prod_i (1-\x_i)}.\label{id1}\ee
Therefore we obtain
\vs{2}

{\bf Theorem} \cite[Theorem 3.1]{TW1}: If $p\ne0$ and $r$ is so small that all poles of the $A_\s$ lie outside $\C{r}$, then
\[\P_Y(x_1(t)=x)=p^{N(N-1)/2}\,\int_{{\C{r}^N}} \prod_{i<j}{\x_j-\x_i\ov f(\x_i,\,\x_j)}
\;{1-\x_1\cd\x_N\ov\prod_i (1-\x_i)}\,\prod_i(\x_i^{x-y_i-1}e^{\ep(\x_i)t})\,d^N\x.\]
\vs{1}
Identity (\ref{id1}) is proved by induction on $N$. Call the left side of the identity $\ph_N(\x_1,\ld,\x_N)$ and the right side $\ps_N(\x_1,\ld,\x_N)$, and assume the identity holds for $N-1$. We first sum over all permutations such that $\s(1)=k$, and then sum over $k$. If we observe that the inequality $i<j$ becomes $j\ne i$ when $i=1$, we see that what we get for the left side of (\ref{id1}), using the induction hypothesis, is
\[{1\ov1-\x_1\,\x_2\cd \x_N}\sum_{k=1}^N(-1)^{k+1}\prod_{j\ne k}
f(\x_k,\,\x_j)\cdot\prod_{j\ne k}\x_j\cdot\ps_{N-1}(\x_1,\ldots,\x_{k-1},\x_{k+1},\ldots,\x_N).\]
If we substitute for  $\ps_{N-1}(\x_1,\ldots,\x_{k-1},\x_{k+1},\ldots,\x_N)$ what it is and do some algebra, we find this would equal the right side of (\ref{id1}) if a simpler identity held:
\be\sum_{k=1}^N\,\prod_{j=1}^Nf(\x_k,\,\x_j)\cdot{1\ov
\x_k\,(p-q\x_k)}\,{1\ov \prod_{j\ne k}(\x_j-\x_k)}={p^{N-1}\ov \prod_j\x_j}-p^{N-1}.\label{simpler}\ee
This one is proved by considering the integral 
\[\int\prod_{j=1}^N(p+qz\x_j-z)\cdot{1\ov
z\,(p-qz)}\cdot{1\ov\prod_{j=1}^N(\x_j-z)}\,dz\]
over a large circle. The integral, and so the sum of the residues at 0, the $\x_k$, and $p/q$, equals zero. This sum is equal to the difference of the two sides of (\ref{simpler}).

\bc{\bf 6. The general particle}\ec

The probability $\P_Y(x_m(t)=x)$ is the sum of $\P(X;t)$ over all $X$ for which $x_m=x$, thus over all $x_1,\ld,x_{m-1}$ statisfying $-\iy<x_1<\cd <x_{m-1}<x$, and all $x_{m+1},\ld,x_{N}$ satisfying $x<x_{m+1}<\cd<x_N<\iy$. The latter we can do, as in the last section, since $r<1$. Eventually we shall expand the $\xs{i}$-contours when $i<m$ to $\C{R}$ with $R>1$ so that we can sum over these $x_i$.

First take a partition $(S_-,\,S_+)$ of $[1,\,N]$ with $|S_-|=m-1$ and sum over all those $\s$ for which $\s([1,\,m-1])=S_-$ and $\s([m,\,N])=S_+$. (At the end we will sum over these partitions.) Set $\s_-=\s|_{[1,\,m-1]},\ \s_+=\s|_{[m,\,N]}$. Then $\s_-$ may be associated in an obvious way with a permutation in $\S_{m-1}$ and  $\s_+$ with a permutation in $\S_{N-m}$. In particular, $\sg\s_\pm$ make sense, and counting inversions shows that
\[\sg\s=(-1)^{\k(S_-,\,S_+)}\,\sg\s_-\,\sg\s_+,\]
where we define in general\footnote{In the cited papers we used the notation $\s(U,\,V)$.}
\[\k(U,\,V)=\#\{(i,\,j):i\in U,\ j\in V,\ i\ge j\}.\]

When we write
\be\prod_i\xs{i}^{x_i}=\prod_{i<m}\xs{i}^{x_i}\ \prod_{i\ge m}\xs{i}^{x_i},\label{prodrep}\ee
the $\s$ in the first product on the right may be replaced by $\s_-$ and the $\s$ in the second product on the right may be replaced by $\s_+$.

Similarly, may rewrite (\ref{A}) as
\be{\ds{\prod_{i<j<m}f(\x_{\s(i)},\,\x_{\s(j)})}\,
\ds{\prod_{i\in S_-,\,j\in S_+}f(\x_i,\,\x_j)}\,
\ds{\prod_{m\le i<j}f(\x_{\s(i)},\,\x_{\s(j)})}\ov
\ds{\prod_{i<j}f(\x_i,\,\x_j)}},\label{A1}\ee
and the $\s$ in the first product may be replaced by $\s_-$ and the $\s$ in the last product may be replaced by $\s_+$

When we sum $\prod_{i\ge m}\xs{i}^{x_i}$ over the $x_i$ with $i\ge m$ (which we may do since $r<1$) we obtain 
\[{\x_{\s(m+1)}\x_{\s(m+2)}^2\cd\x_{\s(N)}^{N-1}\ov (1-\x_{\s(m+1)}\x_{\s(m+2)}\cd\x_{\s(N)}) (1-\x_{\s(m+1)}\cd\x_{\s(N)})\cd(1-\x_{\s(N)})}\;\prod_{i\ge m}\xs{i}^{x}.\]
We then multiply by the last factor in the numerator in (\ref{A1}) times $\sg\s_+$ and sum over $\s_+$. The result is, by (\ref{id1}), 
\be p^{(N-m)(N-m+1)/2}\, {\Big(1-\ds{\prod_{i\in S_+}\x_i}\Big)\,
{\ds\prod_{{i<j\atop i,\,j\in S_+}} (\x_j-\x_i)}\ov \ds{\prod_{i\in S_+} (1-\x_i)}}\,\prod_{i\in S_+}\x_i^x.\label{S+sum}\ee

What remains from (\ref{prodrep}) and (\ref{A1}) is
\be{\ds{\prod_{i<j<m}f(\x_{\s(i)},\,\x_{\s(j)})}\,
\ds{\prod_{i\in S_-,\,j\in S_+}f(\x_i,\,\x_j)}\ov
\ds{\prod_{i<j}f(\x_i,\,\x_j)}}\;\prod_{i<m}\xs{i}^{x_i}.\label{A2}\ee

The next step is to expand all contours to $\C{R}$ with $R>1$, so that we can sum over the $x_i$ with $i<m$. When we expand the contours we encounter poles from the denominators in (\ref{A2}) and (\ref{S+sum}), and integrating the residues would give lower-dimensional integrals. In fact the lower-dimensional integrals from   (\ref{A2}) are zero and the lower-dimensional integrals from (\ref{S+sum}) are of the same type but with fewer variables. Let us see why this is so. We assume $p,\,q\ne0$.

We first consider the poles coming from the denominator in (\ref{A2}) and begin by expanding the $\x_N$ contour to a circle $\C{R}$ with $R$ very large. From the denominator in (\ref{A}) we encounter poles at
\[\x_N={\x_i-p\ov q\x_i}\]
with $i<N$. As in the proof of Lemma 1 we find that the residue at this pole is analytic for $\x_i$ inside $\C{r}$, so the integral with respect to $\x_i$ equals zero. (An important point is that the term $\ep(\x_i)+\ep(\x_N)$, which appears in the exponential in the integrand, becomes analytic at $\x_i=0$ after the substitution.) After expanding the $\x_N$ contour to $\C{R}$ we expand the $\x_{N-1}$-contour. Now from the denominator in (\ref{A2}) we have poles at
\[\x_{N-1}={\x_i-p\ov q\x_i}\]
with $i<N-1$. As before, the integral with respect to $\x_i$ of the residue is equal to zero. There is also the pole at
\[\x_{N-1}={p\ov 1-q\,\x_N}.\]
But $\x_N\in\C{R}$, and if $R$ is chosen large enough this pole is inside $\C{r}$ and so is not crossed in the expansion.

Continuing, we find that when we expand all the contours the poles of $A_\s$ do not contribute. But the poles of (\ref{S+sum}) do contribute, and we get a sum of lower-dimensional integrals, one for each subset  $S_+'\subset S_+$. These are are minus the integrals of the residues at the $\x_k=1$ with $k\in S_+\bs S_+'$. If we use $f(\x_k,\,1)=p\,(1-\x_k)$ and $f(1,\x_k)=q\,(\x_k-1)$ we find that this residue is a constant involving powers of $p$ and $q$ times the same integrand we had before except that there are no terms involving the $\x_k$ with $k\in S_+\bs S_+'$. 

Once all contours are $\C{R}$ (here $R>1$ should be so large that all poles of the $A_\s$ lie inside $\C{R}$) we may sum over all
$x_1,\ld,x_{m-1}$ statisfying $-\iy<x_1<\cd <x_{m-1}<x$. The result of this sum is that the product $\prod_{i<m}\xs{i}^{x_i}$ in (\ref{A2}) is replaced by
\[{1\ov(\xs{1}-1)(\xs{1}\xs{2}-1)\cd(\xs{1}\cd\xs{m-1}-1)}\,\prod_{i\in S_-}\x_i^{x}.\]
Now we are to  multiply by $\sg\s_-$ and sum over all $\s_-$. An identity\footnote{Proved by interchanging $p$ and $q$ in (\ref{id1}) and letting $\x_i\to1/\x_{m-i}$.} analogous to (\ref{id1}) tells us that the sum equals
\[q^{(m-1)(m-2)/2}\, {
{\ds\prod_{{i<j\atop i,\,j\in S_-}} (\x_j-\x_i)}\ov \ds{\prod_{i\in S_-} (1-\x_i)}}\,\prod_{i\in S_-}\x_i^{x}.\]
We now put everything together. We use the notations
\[\t=p/q,\ \ \ \k(U)=\k(U,\,\Z^+)=\textrm{sum of the elements of}  \ U.\]
The result (a special case of \cite[Theorem 3]{TW4}) is that  when $q\ne0$,
\[\P_Y(x_m(t)=x)
=\sum_{S_-,\,S_+'}(-1)^{m-1}\ \t^{\k(S_-\cup S_+')-m(m-1)/2-mk}\,q^{k(k-1)/2+(m-1)(m-2)/2}\]
\[\times\int_{\C{R}}\cd\int_{\C{R}}
{\prod\limits_{i\in S_-,\;j\in S_+'}f(\x_i,\,\x_j)\ 
\ov\prod\limits_{i<j}f(\x_i,\,\x_j)}\ {\Big(1-\ds{\prod_{j\in S_+'}\;\x_j}\Big)\ov\ds{\prod_{i\in S_-,\ j\in S_+'}(\x_j-\x_i)}} 
 {\ds{\prod_{i<j}(\x_j-\x_i)}\ov \ds{\prod_{i}(1-\x_j)}}\ 
 \prod_i(\x_i^{x-y_i-1}\,e^{\ep(\x_i)\,t})\,\prod_i d\x_{i}.\]
Here $k=|S_+'|$ and the sum runs over all disjoint sets $S_-$ and $S_+'$ with $|S_-|=m-1$; indices not otherwise specified run over $S_-\cup S_+'$.

There is one more step. We take a fixed set $S$ and first sum over all partitions $(S_-,\,S_+')$ of $S$ with $|S_-|=m-1$. (At the end we will sum over all these $S$). The only terms that involve $S_-$ and $S_+'$ individually combine as
\[\prod_{i\in S_-,\;j\in S_-^c}{f(\x_i,\,\x_j)\ov \x_j-\x_i}\,
\Big(1-\prod_{j\in S_-^c}\;\x_j \Big),\]
where $S_-^c$ denotes the complement of $S_-$ in $S$. Identity (1.9) of \cite{TW1} (with slightly different notation) is
\be\sum_{{|S_-|=m-1\atop S_-\subset S}}\ \prod_{i\in S_-,\;j\in S_-^c}{f(\x_i,\,\x_j)\ov \x_j-\x_i}\ 
\Big(1-\prod_{j\in S_-^c}\;\x_j \Big)=q^{m-1}\,\br{|S|-1}{m-1}_\t\Big(1-\prod_{i\in S}\x_i\Big),\label{id2}\ee
where the $\t$-{\it binomial coefficient} $\br{N}{n}_\t$ is defined by
\[\br{N}{n}_\t={(1-\t^N)\,(1-\t^{N-1})\cdots (1-\t^{N-n+1})\ov (1-\t)\,(1-\t^2)\cdots (1-\t^n)}.\]
Hence, after some algebra, the result becomes
\vs{2}

{\bf Theorem} \cite[Theorem 5.2]{TW1}: We have when $q\ne0$,
\[\P_Y(x_m(t)=x_m)=(-1)^{m-1}\,\t^{m(m-1)/2}\,\sum_{|S|\ge m}\,\t^{\k(S)-mk}\,q^{k(k-1)/2}\,\br{k-1}{m-1}_\t\]
\be\times \int_{\C{R}^k}\  
\prod_{i<j}{\x_j-\x_i\ov f(\x_i,\,\x_j)}\ {1-\prod_i\x_i\ov 
\prod_i(1-\x_i)}\
\prod_i(\x_{i}^{x-y_i-1}\,e^{\ep(\x_i)\,t})\,d^k\x,\label{mprob}\ee
where now $k=|S|$ and all indices in the integrand run over $S$. 
\vs{1}

Identity (\ref{id2}) depends on a simpler identity,
\be\sum_{{|S_-|=m-1\atop S_-\subset S}}\ \prod_{i\in S_-,\;j\in S_-^c}{f(\x_i,\,\x_j)\ov \x_j-\x_i}=\br{|S|}{m-1}_\t.\label{id3}\ee
This is proved by induction $|S|$. We first observe that the left side is a polynomial in the $\x_i$. The reason\footnote{We learned this argument from Anne Schilling.} is that it is symmetric in the $\x_i$, and if we multiply it by the Vandermonde $\prod_{i<j}(\x_i-\x_j)$ we obtain an antisymmetric polynomial which is, therefore, a polynomial times the Vandermonde. Since the left side it is bounded as each $\x_i\to\iy$ it is a constant. Using the induction hypothesis and a recursion formula for the $\t$-binomial coefficients we see by setting some $\x_i=1$ that the two sides of the identity agree.

To prove (\ref{id2}) by induction, we see now that the left side is a polynomial of degree at most one in each $\x_i$ and, using the induction hypothesis and the recursion formula for the $\t$-binomial coefficients, that the two sides agree when any $\x_i=1$. Therefore the difference is of the form $c\,\prod_i(\x_i-1)$. To show that $c=0$ we use (\ref{id3}) to see that after dividing by some $\x_i$ the two sides of (\ref{id2}) have the same limit as $\x_i\to\iy$. 

\bc{\bf III. Fredholm Determinant Representation for Step Initial Condition}\ec

\bc{\bf 1. Series representation}\ec

Until now we assumed a system of finitely many particles. Because we can take arbitrarily large $R$ in (\ref{mprob}) the result extends to initial configurations
\[y_1<y_2<\ld\to+\iy,\]
where the sum is taken over all finite sets $S\subset\Z^+$.

For step initial configuration, where $Y=\Z^+$ and $y_i=i$, we may sum over all sets $S$ with $|S|=k$ and so obtain instead a sum over $k\ge m$. Before that, instead of indexing the variables in the integrand by $S$ we index them by $\{1,\ld,k\}$, so that we can sum under the integral signs for all $S$ with $|S|=k$. If $S=\{s_1,\ld,s_k\}$ with $s_1<\cd<s_k$ then after renaming the variables the factor $\prod_i \x_i^{-y_i}$ in the integrand becomes $\prod_i \x_i^{-s_i}$ and $\t^{\k(S)}$ becomes $\prod_i \t^{s_i}$. These are the only terms in (\ref{mprob}) that involve the individual $s_i$. Summing the product of these two over all $s_i$ satisfying $1\le s_1<\cd<s_k<\iy$ gives, when $R>\t$, 
\be{1\ov \Big((\x_1/\t)
(\x_2/\t)\cd(\x_k/\t)-1\Big)\,\Big
((\x_2/\t)\cd(\x_k/\t)-1\Big)\cd\Big((\x_k/\t)-1\Big)}.\label{ksum}\ee
The factor $\prod_{i<j}f(\x_i,\,\x_j)\inv$ in the integrand may be written
\[{\ds{\prod_{i>j}f(\x_i,\,\x_j)}\ov\ds{\prod_{i\ne j}f(\x_i,\,\x_j)}}.\]
If we multiply (\ref{ksum}) by the numerator here the rest of the integrand is antisymmetric in the $\x_i$. Thus the integral is unchanged if this product is antisymmetrized. We make the substitution $\x_i\to \t/\x_i$, use identity (\ref{id1}) with $p$ and $q$ interchanged, and find that the antisymmetrization is
\[{1\ov k!}\,p^{k(k+1)/2}\,{\ds{\prod_{i>j}(\x_j-\x_i)}\ov\ds{\prod_i(q\x_i-p)}}.\]
Thus we obtain, 
\vs{2}

{\bf Theorem} \cite[Corollary to Th. 5.2]{TW1}: For step initial condition we have when $q\ne 0$,

\[\P(x_m(t)\le x)=(-1)^{m}\,
\sum_{k\ge m}{1\ov k!}\,\br{k-1}{k-m}_\t\,\t^{m(m-1)/2-mk+k/2}\;
(pq)^{k^2/2}\]
\[\times\int_{\C{R}}\cd\int_{\C{R}}\prod_{i\ne j}{\x_j-\x_i\ov f(\x_i,\,\x_j)}\;\prod_i{\x_i^{x}\,e^{\ep(\x_i)t}\ov(1-\x_i)\,(q\x_i-p)}\,d\x_1\cd d\x_k.\]
\vs{1}
Notice that on the left side we have $\P(x_m(t)\le x)$ rather than $\P(x_m(t)=x)$ and on the right side the sign is different and the factor $1-\prod_i\xi_i$ is gone. This is the result of summing the formula for $\P(x_m(t)=x)$ from $-\iy$ to $x$.

For TASEP with $p=0$ only the term $k=m$ is nonzero, the multiple integral is an $m\times m$ Toeplitz determinant, and we get 

\[\P(x_m(t)\le x)=\det\(\int_{\C{R}}\x^{i-j+x-1}\,(\x-1)^{-m}\,e^{(\x-1)t}\,d\x\).\]
This was obtained by R\'akos and Sch\"utz \cite{RS} who showed it was equivalent to Johansson's result mentioned in the introduction.

\bc{\bf 2. Fredholm determinant representation}\ec

If we make the change of variables
\[\x_i={1-\t\e _i\ov 1-\e_i},\]
then
\[\prod_{i\neq j}{\x_j-\x_i\ov p +q \x_i\x_j -\x_i}=(1+\t)^{k(k-1)}\,\prod_{i\neq j}{\e_i-\e_j\ov \t\e_i-\e_j}.\]
The right side can be represented in terms of the Cauchy determinant
\[\det\({1\ov \t\e_i-\e_j}\)=\t^{k(k-1)/2}\,{\prod_{i\ne j}(\e_i-\e_j)\ov \prod_{i,j}(\t\e_i-\e_j)}.\]
Going back to the $\x_i$ gives the identity
\[\prod_{i\neq j}{\x_j-\x_i\ov p +q \x_i\x_j -\x_i}=(-1)^k\,(pq)^{-k(k-1)/2}\,\prod_i (1-\x_i)(q\x_i-p)\,\cdot\, 
\det\left({1\ov p + q \x_i\x_j - \x_i}\right)_{1\le i,j\le k}.\]
The theorem becomes
\[\P(x_m(t)\le x)=(-1)^{m}\,\t^{m(m-1)/2}
\sum_{k\ge m}\br{k-1}{k-m}_\t\,\t^{(1-m)k}\]
\[\times{(-1)^k\ov k!}\,\int_{\C{R}}\cd\int_{\C{R}}\det (K(\x_i,\,\x_j))_{1\le i,j\le k}\;d\x_1\cd d\x_k,\]
where
\[K(\x,\,\x')=q\,{\x^x\,e^{\ep(\x)t}\ov p+q\x\x'-\x}.\]

Denote by $K$ the operator acting on functions on $\C{R}$ by 
\[Kf(\x)=\int_{\C{R}}K(\x,\x')\,f(\x')\,d\x'.\]
The Fredholm expansion is
\[\det(I-\la K)=\sum_{k=0}^\iy {(-\la)^k\ov k!}\,\int_{\C{R}}\cd\int_{\C{R}}\det K(\x_i,\,\x_j)_{1\le i,j\le k}\;d\x_1\cd d\x_k,\]
which gives
\[{(-1)^k\ov k!}\,\int_{\C{R}}\cd\int_{\C{R}}\det K(\x_i,\,\x_j)_{1\le i,j\le k}\;d\x_1\cd d\x_k=\int {\det(I-\la K)\ov\la^{k+1}}\,d\la,\] 
where we take any contour enclosing $\la=0$. Thus,
\[\P(x_m(t)\le x)=(-1)^{m}\,\t^{m(m-1)/2}
\sum_{k\ge m}\br{k-1}{k-m}_\t\,\t^{(1-m)k}\,\int {\det(I-\la K)\ov\la^{k+1}}\,d\la.\]
If the contour is $\C{\rho}$ with $\rho>\t^{1-m}$ then we can interchange the sum and integral and use the $\t$-binomial theorem 
\[\sum_{k\ge m} \left[{k-1\atop k-m}\right]_\t z^{k}=\prod_{j=1}^{m}{z\ov  1- \t^{m-j} z}\]
with $z=\t^{1-m}\,\la\inv$.
We obtain,
\vs{2}

{\bf Theorem} \cite[Formula (1)]{TW2}: We have when $q\ne0$,
\be \P\left(x_m(t)\le x\right)=\int {\det(I-\la K)\ov\prod_{k=0}^{m-1}(1-\la\,\t^k)}\, {d\la\ov \la},\label{Fredrep}\ee
where the contour of integration encloses all the singularities of the integrand.
\vs{1}

We can evaluate the integral by residues, getting a finite sum of determinants. When $m=1$ we obtain
\[\P\left(x_1(t)>x\right)=\det(I- K),\]

\begin{center}{\bf IV. Asymptotics}\end{center}

\bc{\bf 1. Statements of the rsults}\ec

If $q>p$ and $t\to\iy$, we expect $x_m(t)$ to be large and negative. We obtain three asymptotic results for $\P\left(x_m(t)\le x\right)$ as $t\to\iy$. Recall the definition $\g=q-p$.
\vs{2}

{\bf Theorem 1} \cite[Theorem 1]{TW3}: Let $m$ and $x$ be fixed with $x<m$. Then as $t\to\iy$
\[\P\left(x_m(t/\g)>x\right)\sim \prod_{k=1}^\iy (1-\t^k)\,{t^{2m-x-2}\,e^{-t}\ov (m-1)!\,(m-x-1)!}.\]
\vs{2}

{\bf Theorem 2} \cite[Theorem 2]{TW3}: For fixed $m$ we have 
\[\lim_{t\to\iy}\P\(x_m(t/\g)\le -t-\g^{1/2}\,s\,t^{1/2}\)=\int {\det(I-\la \hat{K}_s)\ov\prod_{k=0}^{m-1}(1-\la\,\t^k)} {d\la\ov \la},\]
where $\hat{K}_s$ is the operator on $L^2(s,\iy)$ with kernel
\[\hat{K}(z,\,z')={q\ov\sqrt{2\pi}}\,\,e^{-(p^2+q^2)\,(z^2+{z'}^2)/4+pq\,zz'}.\]

For the third result we recall that 
\[F_2(s)=\det\,\big(I-K_{{\rm Airy}}\,\ch_{(s,\iy)}\big),\]
where 
\[K_{\rm Airy}(x,\,y)=\int_0^\iy\A(x+z)\,\A(y+z)\,dz.\]
For $\s\in(0,1)$ we set
\be c_1(\s) =-1+2\sqrt{\s},\ \ \ c_2(\s)=\s^{-1/6}\,(1-\sqrt{\s})^{2/3}.\label{sigma}\ee
\vs{2}

{\bf Theorem 3} \cite[Theorem 3]{TW3}: We have \[\lim_{t\to\iy}\P\({x_{[\s t]}(t/\g)\le c_1(\s)\,t+ c_2(\s)\,s\,t^{1/3}}\)=F_2(s)\]
uniformly for $\s$ in compact subsets of $(0,\,1)$.
\vs{2}
For TASEP ($p=0$) the probabilities are $m\times m$ determinants. For $m$ and $x$ fixed the asymptotics of the determinant are easily found and agree with Theorem 1.

A special case of Theorem 2 is
\[\lim_{t\to\iy}\P\(x_1(t/\g)>-t-\g^{1/2}\,s\,t^{1/2}\)=\det(I-\hat{K}_{s}).\]
This is a family of distribition functions parameterized by $p\in[0,\,1)$. When $p=0$ it is a normal distribution and the probability on the left is the probability for a free particle.

Theorem 3 when $p=0$ gives the asymptotics for TASEP obtained by Johansson~\cite{J}. A consequence of Theorem 3 is that for fixed $s\in(0,\,1)$
\[\lim_{t\to \iy}{x_{[\s t]}(t/\g)\ov t}=c_1(\s).\]
in probability. In fact Liggett \cite{L} showed that this holds almost surely.
\pagebreak

\begin{center}{\bf 2. Preliminaries}\end{center}

A natural approach to the asymptotics is to look for a limiting operator $K_\iy$ such that $\det(I-\la K)\to\det(I-\la K_\iy)$ as $t\to\iy$. Once one has guessed $K_\iy$ there are two possible approaches:

(i) Show that $K\to K_\iy$ in trace norm.

(ii) Show that $\tr\, K^n\to\tr\, K_\iy^n$ for each $n\in\Z^+$ and that $K$ is bounded in trace norm (or even Hilbert-Schmidt norm) as $t\to\iy$. This suffices because of the general formula
\be\log\det(I-\la L)=-\sum_{n=1}^\iy {\la^n\ov n}\,\tr\,L^n,\label{logdet}\ee
which holds for sufficiently small $\la$. 
 
Both approaches will be used eventually. The operators $K$ on $\C{R}$ have exponentially large norms as $t\to\iy$, and we will replace them by operators with the same Fredholm determinants that are better-behaved. This will be possible because of lemmas on
stability of Fredholm determinants.
\vs{1}
{\bf Lemma 1}. If $s\to\G_s$ is a deformation of closed curves and 
$L(\e,\,\e')$ is analytic for $\e,\,\e'\in \G_s$ for all $s$, then the Fredholm determinant of $L$ acting on $\G_s$ is independent of $s$.
\vs{2}
{\bf Lemma 2}. If $L_1(\e,\,\e')$ and $L_2(\e,\,\e')$ are two kernels acting on a simple closed curve $\G$, such that $L_1(\e,\,\e')$ extends analytically to $\e$ inside $\G$ {\bf or} to $\e'$ inside $\G$, and $L_2(\e,\,\e')$ extends analytically to $\e$ inside $\G$ {\bf and} to $\e'$ inside $\G$, then the Fredholm determinants of $L_1(\e,\,\e')+L_2(\e,\,\e')$ and $L_1(\e,\,\e')$ are equal.
\vs{1}
The proofs use the fact that $\det(I-\la L)$ is determined by the traces $\tr\, L^n$. For Lemma 1 we use
\[\tr\,L^n=\int_{\G_s}\cdots\int_{\G_s}L(\e_1,\,\e_2)\cdots L(\e_{n-1},\,\e_n)\,L(\e_n,\,\e_1)\,d\e_1\cdots d\e_n.\]
By analyticity the integral is invariant under the deformation.
For Lemma 2, we have to show
\[\tr\,(L_1+L_2)^n=\tr\,L_1^n.\]
If, say, $L_1(\e,\,\e')$ extends analytically to $\e'$ inside $\G$, then 
\[L_1\,L_2\,(\e,\,\e'')=\int_\G L_1(\e,\,\e')\,L_2(\e',\,\e'')\,d\e'=0.\]
Similarly $L_2^2=0$. Also, $\tr\,L_2=0$, so $\tr\,(L_1+L_2)=\tr\,L_1$. 
Since $L_1L_2=L_2^2=0$, we have for $n>1$ 
\[(L_1+L_2)^n=L_1^n+L_2\,L_1^{n-1}.\]
Since 
\[\tr\,L_2\,L_1^{n-1} =\tr\,L_1^{n-1}\,L_2=0,\]
we have $\tr\,(L_1+L_2)^n=\tr\,L_1^n$.
\vs{2}
\bc{\bf 3. Another operator}\ec

We introduce the notation
\[\ph(\e)=\({1-\t \e\ov1-\e}\)^x\,e^{\left[{1\ov 1-\e}-
{1\ov 1-\t \e}\right]\,t}.\]
In $K(\x,\,\x')$ we make the substitutions
\[\x={1-\t\e\ov1-\e},\ \ \ \x'={1-\t\e'\ov1-\e'}, \ \ \ t\to t/\g\]
and obtain the kernel\footnote{This is the kernel $(d\x/d\e)^{1/2}(d\x'/d\e')^{1/2}\,K(\x(\e),\,\x'(\e'))$.} 
\[{\ph(\e')\ov \e'-\t\e}=K_2(\e,\,\e')\]
acting on \c, a little circle about $\e=1$ described clockwise, which has the same Fredholm determinant. We denote this by $K_2$ because there is an equally important kernel 
\[{\ph(\t\e)\ov \e'-\t\e}=K_1(\e,\,\e').\]

The kernel $K_1(\e,\,\e')$ extends analytically to $\e$ and $\e'$ inside \c\  while $K_2(\e,\,\e')$ extends analytically to $\e$ inside \c. Hence by Lemma~2 the determinant of $K_2$ equals the determinant of $K_2-K_1$.

Next we apply Lemma~1 to the kernel 
\[K_1(\e,\,\e')-K_2(\e,\e')={\ph(\t\e)-\ph(\e')\ov \e'-\t\e},\]
with $\G_0=-$\c\ and $\G_1=\C{\r_\e}$ with $1<\r_\e<\t\inv$. (Recall that \c\ was described clockwise.) Since the numerator vanishes when the denominator does, the only singularities of the kernel are at $\e,\,\e'=1,\,\t\inv$, neither of which is passed in a deformation $\G_s,\ s\in[0,\,1]$. Therefore the operator $K$ acting on  $\C{R}$ may be replaced by $K_1-K_2$ acting on $\C{\r_\e}$.
\pagebreak

\bc{\bf 4. Another Fredholm determinant representation}\ec 

The function $\ph(\t\e)$ is analytic on $s\C{\r}$ when $0<s\le1$. The denominator $\e'-\t\e$ in $K_1$ is nonzero for $\e,\,\e'\in s\C{\r}$ for all such $s$. Therefore by Lemma~1 the Fredholm determinant of $K_1$ on $\C{\r}$ is the same as on $s\C{\r}$. This in turn is the same as the Fredholm determinant of
\be s\,K_1(s\e,\,s\e')={\ph(s\t\e)\ov \e'-\t\e}\label{Ks}\ee
on $\C{\r}$. As $s\to0$ this converges in trace norm to the kernel
 \[K_0(\e,\,\e')={1\ov \e'-\t\e}\]
on $\C{\r}$. Therefore the Fredholm determinant of $K_1$ equals the Fredholm determinant of $K_0$.

The kernel of $K_0^2$ equals
\[ K_0^2(\e,\,\e')=\int_{\C{\r}}{d\z\ov (\z-\t \e)\,(\e'-\t \z)}=
{1\ov \e'-\t^2\,\e},\]
because $\t\e$ is inside $\C{\r}$ and $\t\inv\e'$ outside when  $\e,\,\e'\in\C{\r}$. In particular $\tr\,K_0^2=(1-\t^2)\inv$. Generally, we find that $\tr\,K_0^n=(1-\t^2)^{-n}$. Thus by (\ref{logdet}) we have for small $\la$
\[\log\det(I-\la K_0)=-\sum_{n=1}^\iy{\la^n\ov n}{1\ov1-\t^n}=-\sum_{k=0}^\iy\sum_{n=1}^\iy {\t^{nk}\la^n\ov n}=\sum_{k=0}^\iy\log(1-\la\t^k),\]
and so
\[\det(I-\la K_1)=\det(I-\la K_0)=\prod_{k=0}^\iy (1-\la\t^k).\] 

We  factor out $I-\la K_1$ in
\[\P(x_m(t/\g)\le x)=\int {\det(I-\la K)\ov\prod_{k=0}^{m-1}(1-\la\,\t^k)}\,{d\la\ov \la}=\int {\det(I-\la K_1+\la K_2)\ov\prod_{k=0}^{m-1}(1-\la\,\t^k)}\,{d\la\ov \la},\]
(recall the substitution $t\to t/\g$) and obtain 
\be\P(x_m(t/\g)\le x)=\int \prod_{k=m}^{\iy}(1-\la\,\t^k)\,\det(I+\la K_2\,(I+R))\,{d\la\ov \la},\label{probform2}\ee
where $R$ is the resolvent operator $\la K_1\,(I-\la K_1)\inv$.
\pagebreak  

\bc{\bf 5. Theorems 1 and 2}\ec

Consecutive integration shows that for small $\la$ the resolvent kernel has the nice representation
\be R(\e,\,\e';\,\la)=\sum_{n=1}^\iy \la^n {\ph_n(\t\e)\ov \e'-\t^n\e},\label{R}\ee
where 
\[\ph_n(\e)=\ph(\e)\,\ph(\t\e)\cdots\ph(\t^{n-1}\e).\]
For Theorems 1 and 2, whose derivations we shall not explain in detail, we wrote $R=R_1+R_2$ where $R_1$ is analytic everywhere except for the poles at $\la=1,\,\t\inv,\ldots,\t^{-m+1}$ and $R_2$ is analytic for $|\la|<\t^{-m}$. For Theorem~1 the asymptotics comes from the residue of $R_1$ at $\la=\t^{-m+1}$. For Theorem~2 we used approach (ii) described above. In \cite{TW2} a steepest descent computation had shown that $\tr\,K^n\to \tr\,\hat{K}^n$ for all $n$. What was needed to complete the proof was to show that $K_2\,(I+R)$ had bounded Hilbert-Schmidt norm as $t\to\iy$, uniformly for $\la$ in compact sets not containing any of the singularities $\t^{-k}$. We used the representation $R=R_1+R_2$ to show that this was so.

\bc{\bf 6. Theorem 3}\ec
Here $m=\s t$ is large and $1,\,\t\inv,\ld,\t^{-m+1}$ must be inside the contour. If we set $\la=\m\,\t^{-m}$ we can take $\m\in\C{\r}$ with $\r>\t$ fixed, and (\ref{probform2})  becomes
\be\P(x_m(t/\g)\le x)=\int \prod_{k=0}^{\iy}(1-\m\,\t^k)\ \det(I+\m\,\t^{-m} K_2\,(I+R))\,{d\m\ov \m}.\label{probform3}\ee

In (\ref{R}) we use
\[\ph_n(\e)={\phy(\e)\ov\phy(\t^n\,\e)},\]
where 
\[\ph_\iy(\e)=\lim_{n\to\iy}\ph_n(\e)=(1-\e)^{-x}\,e^{{\e\ov 1-\e}t}.\]
The Cauchy integral representation of $\phy(\t^n\,\e)\inv$, and some manipulation of series and integrals, give
\[K_2(I+R)\,(\e,\,\e')=-\int_{|\z|>\r_\e}{\ph(\z)\ov (\z-\t\e)\,(\e'-\z)}\,d\z\]
\[+\sum_{k=-\iy}^\iy{\t^{k}\ov 1-\t^k\la}
\int_{\C{\r_\z}}{\phy(\z)\ov\z-\t\e}\,\z^k\,d\z\,\int_{\C{\r_u}}{1\ov\phy(u)\,(\e'-u/\t)}\;{du\ov u^{k+1}},\]
where the radii of the contours in the series satisfy
\[\r_\z\in(1,\,\r_\e),\ \ \ \r_u\in (\t\,\r_\z,\,\t\,\r_\e).\]
 
The first operator on the right side  is analytic for $|\e|,\,|\e'|\le \r_\e$ and the infinite sum is analytic for $|\e|\le \r_\e$. It follows by Lemma~2 that the Fredholm determinant of the sum of the two, i.e., of $K_2(I+R)$, equals the Fredholm determinant of the infinite sum.

If we set
\[f(\m,\,z)=\sum_{k=-\iy}^\iy{\t^{k}\ov 1-\t^{k}\m}\,z^k,\]
then since $\la=\m\,\t^{-m}$,  
\[\sum_{k=-\iy}^\iy{\t^{k}\ov 1-\t^k\la}\({\z\ov u}\)^k=\t^m\,\({\z\ov u}\)^m\,f(\m,\z/u),\]
and so the infinite sum may be written 
\[\t^m\,\int_{\C{\r_u}}\int_{\C{\r_\z}}{\phy(\z)\ov\phy(u)}\,\({\z\ov u}\)^m\,{f(\m,\z/u)\ov(\z-\t\e)\,(\e'-u/\t)}\;d\z\;{du\ov u}.\]

The substitutions $\e\to\e/\t,\ \e'\to\e'/\t$ replace this by the kernel 
\[\t^m\,\int_{\C{\r_u}}\int_{\C{\r_\z}}{\phy(\z)\ov\phy(u)}\,\({\z\ov u}\)^m\,{f(\m,\z/u)\ov(\z-\e)\,(\e'-u)}\;d\z\;{du\ov u},\]
where now the operator acts on $\C{\r_\e}$ with $\r_\e\in (\t,\,1)$ and in the integral
\[\r_\z\in(1,\,\t\inv\,\r_\e),\ \ \ \r_u\in(\t\r_\z,\,\r_\e).\]

If we expand the $u$-integral so that $\r_\e<|u|<1$ on the new contour we pass the pole at $u=\e'$, which gives the contribution
\be\t^m\,\int_{\C{\r_\z}}{\phy(\z)\ov\phy(\e')}\,{\z^m\ov (\e')^{m+1}}\,{f(\m,\z/\e')\ov\z-\e}\;d\z.\label{contribution}\ee
The new double integral is analytic for $|\e|,\,|\e'|\le \r_\e$ and (\ref{contribution}) is analytic for $|\e|\le \r_\e$. Therefore by Lemma 2 the Fredholm determinant is the same as that of (\ref{contribution}). 

We have shown that if we define
\be J(\e,\,\e')=\int_{\C{\r_\z}}{\phy(\z)\ov\phy(\e')}\,{\z^m\ov (\e')^{m+1}}\,{f(\m,\z/\e')\ov\z-\e}\;d\z,\label{J}\ee
where $\r_\z\in(1,\,\t\inv\,\r_\e)$, then (\ref{probform3}) becomes
\be\P(x_m(t/\g)\le x)=\int \prod_{k=0}^{\iy}(1-\m\,\t^k)\ \det(I+\m\,J)\,{d\m\ov \m}.\label{probform4}\ee
This representation, in which the parameter $m$ is in the operator, makes an asymptotic analysis possible.

By Lemma 1 the contours $\C{\r_\e}$ (the home of the functions on which $J$ acts) and $\C{\r_\z}$ (in the integral defining $J$) may be simultaneously deformed if during the deformation we do not pass a singularity of the integrand.

We apply steepest descent, and so look for the saddle points for
$\ph(\z)\,\z^m$ when $m\sim\s t$ and $x\sim c\,t$. In general there are two saddle points. When $c$ equals $c_1(\s)$, given in (\ref{sigma}), they coincide at
\[\x=-\sqrt\s/(1-\sqrt{\s}).\]
Both contours may be deformed to pass through the saddle point, the neighborhood of which gives the main contributions. If $x=c_1(\s)\,t+ c_2(\s)\,s\,t^{1/3}$ precisely, and we make the variable changes 
\[\e\to \x+t^{-1/3}\,c_3\,\e,\ \ \ \e'\to \x+t^{-1/3}\,c_3\,\e',\ \ \ \z\to \x+t^{-1/3}\,c_3\,\z\]
for a certain constant $c_3$, then the rescaled kernel $\m\,J(\m,\,\m')$ has limit
\[\int_{\G_\z} {e^{-\z^3/3+s\z+(\e')^3/3-s\e'}\ov(\z-\e)\,(\e'-\z)}\,d\z.\]
(The constants $c_2(\s)$ and $c_3$ come from a third derivative at the saddle point.) Here $\G_\z$ consists of the the rays from $0$ to $\iy\,e^{\pm2\pi i/3}$. The limiting operator acts on functions on the contour $\G_\e$ consisting of the the rays from $0$ to $\iy\,e^{\pm \pi i/3}$. 

Using the general identity $\det(I-AB)=\det(I-BA)$ we replace this by the kernel
\[
\int_{\G_\z}\int_{\G_\e}{e^{-\z^3/3+\e^3/3+y\z-x\e}\ov\z-\e}\,d\e\,d\z
=-\KA(x,\,y),\]
acting on $L^2(s,\,\iy)$, where
\[\KA(x,\,y)=\int_0^\iy \A(z+x)\,\A(z+y)\,dz.\footnote{The reason the double integral equals $-\KA(x,y)$ is that applying the operator $\partial/\partial x+\partial/\partial y$ to the two kernels
gives the same result, $\A(x)\,\A(y)$, so they differ by a function of $x-y$. Since both kernels go to zero as $x$ and $y$ go to $+\iy$ independently this function must be zero.}\]

Hence
\[\det(I+\m\, J)\to\det\(I-\KA\,\ch_{(s,\,\iy)}\)=F_2(s)\]
for all $\m$, and it follows that the integral in (\ref{probform4}) has the limit $F_2(s)$.

\begin{center}{\bf Acknowledgment}\end{center}

The authors were supported by the National Science Foundation through grants DMS-0906387 (first author) and DMS-0854934 (second author).

\end{document}